\documentclass[11pt]{article}
\begin{document}
\tolerance=10000


\hyphenpenalty=2000
\hyphenation{visco-elastic visco-elasticity}
\setcounter{page}{1}
\thispagestyle{empty}




\font\note=cmr10 at 10 truept  
\font\note=cmr8  

\newcommand{\FTS}[2]{\frac{{\textstyle #1}}{{\textstyle #2}}}
\newcommand{\nat}{I\!\!N}
\newcommand{\NN}{I\!\!N}
\newcommand{\ka}{I\!\!K}
\newcommand{\rgr}{{\rm grad}}
\newcommand{\ce}{I\!\!\!\!C}
\newcommand{\CC}{I\!\!\!\!C}
\newcommand{\re}{I\!\!R}
\newcommand{\RR}{I\!\!R}
 \newcommand{\intl}{\int\limits}
\newcommand{\suml}{\sum\limits}
\setcounter{page}{1}
\thispagestyle{empty}
\font\bfs=cmbx10 scaled \magstep2
\def\eg{{\it e.g.}\ } \def\ie{{\it i.e.}\ }
\def\sg{\hbox{sign}\,}
\def\sgn{\hbox{sign}\,}
\def\sign{\hbox{sign}\,}
\def\e{\hbox{e}}
\def\exp{\hbox{exp}}
\def\ds{\displaystyle}
\def\dis{\displaystyle}
\def\q{\quad}    \def\qq{\qquad}
\def\lan{\langle}\def\ran{\rangle}
\def\l{\left} \def\r{\right}
\def\lra{\Longleftrightarrow}
\def\arg{\hbox{\rm arg}}
\def\d{\partial}
 \def\dr{\partial r}  \def\dt{\partial t}
\def\dx{\partial x}   \def\dy{\partial y}  \def\dz{\partial z}
\def\rec#1{{1\over{#1}}}
\def\log{\hbox{\rm log}\,}
\def\erf{\hbox{\rm erf}\,}     \def\erfc{\hbox{\rm erfc}\,}
\def\G{{G_{\alpha,\beta}^\theta}}
\def\K{K_{\alpha,\beta}^\theta}
\def\Gxt{\G (x,t)}
\def\Gkt{{\widehat{\G}}  (\kappa,t)}
\def\Gxs{{\widetilde{\G}}  (x,s)}
\def\Gks{{\widehat{\widetilde {\G}}} (\kappa,s)}
\def\FT{{\cal F}\,} 
\def\LT{{\cal L}\,}  
\def\L{{\cal L}} 
\def\F{{\cal F}} 
\def\M{{\cal M}}  
\def\I{{\cal I}}  
\def\pni{\par \noindent}
\def\vsh{\smallskip}
\def\vs{\medskip}
\def\vvs{\bigskip}
\def\vvvs{\bigskip\medskip} 
\def\vsp{\par}
\def\vsn{\vsh\pni}
\def\cen{\centerline}
\def\ra{\item{a)\ }} \def\rb{\item{b)\ }}   \def\rc{\item{c)\ }}
\def\alphak{{\alpha \choose k}}
\def\alphazero{{\alpha \choose 0}}
\def\alphaone{{\alpha \choose 1}}
\def\alphatwo{{\alpha \choose 2}}
\def\alphakk{{\alpha \choose k+1}}

 \markboth
 {\rm \centerline {F. Mainardi}}
 {\rm \centerline{INTEGRAL TRANSFORMS IN FRACTIONAL DIFFUSION}}
$\null$
\cen{{\bf FRACALMO PRE-PRINT} \    {\tt www.fracalmo.org}}
\vsh
\cen{{\bf Integral Transforms and Special Functions,}}
\cen{{\bf  Vol 15, No 6, pp. 477-484 (2004)}}
\vsh
\hrule
   \vskip 1.50truecm

\font\title=cmbx12 scaled\magstep1
\font\bfs=cmbx12 scaled\magstep1


\centerline{{\title APPLICATIONS OF INTEGRAL TRANSFORMS}}
\vvs

\centerline{{\title IN FRACTIONAL DIFFUSION PROCESSES}}

\vvs

\begin{center}

 {FRANCESCO \ MAINARDI}  


\vs

 {\it Dipartimento di Fisica, Universit\`a di Bologna and INFN,} \\
{\it Via Irnerio 46, I-40126 Bologna, Italy} \\
 E-mail: {\tt mainardi@bo.infn.it} $\;$ URL: {\tt www.fracalmo.org}

\end{center}

\vvs

{\baselineskip =11pt    \noindent
The fundamental solution (Green function) for the Cauchy problem
of the space-time fractional diffusion equation
is investigated with respect to its scaling and similarity properties,
starting from its Fourier-Laplace representation.
Then,
 by using the Mellin transform,  a general representation
of the  Green function in terms of
Mellin-Barnes integrals in the complex plane is derived.
This
allows us to obtain a suitable computational form of the Green function
in the space-time domain and to analyse its  probability
interpretation.}


\rm 
\section{INTRODUCTION}

In this paper, we review the Cauchy problem for the
{\it space-time fractional} partial differential equation,
which is obtained from the standard diffusion equation
by replacing the second-order space derivative
with a Riesz-Feller derivative of order $\alpha \in (0,2]$
and skewness $\theta$ ($|\theta|\le\hbox{min}\,\{\alpha ,2-\alpha \}$),
and the first-order time derivative with a Caputo derivative
of order $\beta \in (0,2]\,.$

\catcode`@=11  
\font\fnotefont=cmr10 scaled \magstephalf
\def\fnote#1#2{\let\@sf\empty 
  \ifhmode\edef\@sf{\spacefactor=\the\spacefactor}\/\fi
  {$\,^{#1}$}\@sf\vfootnote{#1}{#2}}
\def\vfootnote#1#2{\insert\footins\bgroup
  \interlinepenalty\interfootnotelinepenalty
  \splittopskip=\ht\strutbox 
  \splitmaxdepth=\dp\strutbox \floatingpenalty=20000
  \ifdim\lastskip<2truept \removelastskip\vskip 2truept\fi
  \leftskip=\parindent \rightskip=0pt \spaceskip=0pt \xspaceskip=0pt
  \baselineskip=11pt\parindent=0pt\fnotefont
  \textindent{$\,^{#1}$}#2 \footstrut\vskip 3truept \egroup}
\def\footstrut{\vbox to\splittopskip{}}
\vfill\eject

The fundamental solution (Green function) for the Cauchy problem is
investigated with respect to its scaling and similarity properties,
starting from its Fourier-Laplace representation.
In the  cases
 $\,\{0<\alpha\le 2\,,\,\beta =1\}$ and
 $\,\{\alpha=2\,,\,0<\beta \le 1 \}\,,$
the fundamental solution is known to be interpreted
as a {\it spatial  probability density function evolving in time},
so we talk of {\it space-fractional diffusion} and
{\it time-fractional diffusion}, respectively.
Then, by using the Mellin transform, we provide a general representation
of the  Green functions in terms of Mellin-Barnes integrals in the complex
plane, which allows us {\it to extend the  probability interpretation
to the ranges} $\{0<\alpha \le 2\,,\, 0<\beta \le 1 \}$
{\it and} $\{1<\beta \le \alpha \le 2\} $.


Furthermore, from this representation it is possible to derive
explicit formulae (convergent series and asymptotic expansions),
which  enable us to plot the spatial probability densities
for different values of the relevant
parameters $\alpha, \theta, \beta\,.$

\section{THE GREEN FUNCTION}

By  replacing in the standard diffusion equation
$$ {\d\over \dt} u(x,t) =  {\d^2\over \dx^2}\,
  u(x,t)\,,
   \q -\infty< x <+\infty\,, \q t \ge 0\,,
   \eqno(1)$$
where $u=u(x,t)$ is the (real) field variable,
the second-order space derivative and the first-order time derivative
by    suitable {\it integro-differential} operators, which can be
interpreted as a space and time derivative of fractional order,
we  obtain a   sort of "generalized diffusion" equation.
Such equation  may be
referred to as the {\it space-time-fractional diffusion} equation
when its fundamental solution (see below) can be interpreted
as a probability density.
We write 
$$
_tD_*^\beta \, u(x,t) \,       = \, _xD_\theta^\alpha \,u(x,t) \,,
\q -\infty< x <+\infty\,, \q t \ge 0\,,
\eqno(2) $$
where  the $\alpha \,,\,\theta\,,\, \beta $ are real parameters
restricted as follows
$$ 0<\alpha\le 2\,,\q |\theta| \le \hbox{min} \{\alpha, 2-\alpha\}\,,
  \quad 0<\beta \le 2\,.\eqno(3)$$
In Eq. (2)
$\,_xD_\theta^\alpha \,$ is
the space-fractional
{\it Riesz-Feller derivative}  of order $\alpha $ and skewness $\theta\,,$
and  $\,_tD_*^\beta\,$  is
 the time-fractional {\it Caputo derivative} of order $\beta \,.$
The definitions of these fractional derivatives
are more easily understood if given
in terms of Fourier transform and Laplace transform, respectively.

\noindent
For the space-fractional {\it Riesz-Feller derivative} we have
$$ \qq \qq \qq \qq  \qq
{\cal F} \l\{\, _xD_\theta^\alpha\, f(x);\kappa \r\} =
  - \psi_\alpha ^\theta(\kappa ) \,
  \, \widehat f(\kappa) \,, \qq \qq \qq \qq \qq \eqno(4)$$
$$
   \psi_\alpha ^\theta(\kappa ) =
|\kappa|^\alpha \, \e^{\ds  i (\sgn \kappa)\theta\pi/2}\,,
\q \kappa \in \RR\,,  $$
 where
$  \widehat f(\kappa)  =
{\cal F} \l\{ f(x);\kappa \r\}
  = \int_{-\infty}^{+\infty} \e^{\,\ds +i\kappa x}\,f(x)\, dx\,.$
In other words the symbol of the pseudo-differential operator
\footnote{
Let us recall that
a generic linear pseudo-differential operator $A$,
acting with respect to the variable $x \in \RR\,,$
is defined through its Fourier representation, namely
  $    \int_{-\infty}^{+\infty}
  \e ^{\, i\kappa x} \,  A \,[ f(x)] \, dx =
 \widehat A(\kappa )\, \widehat f (\kappa )\,,  $
  where
$\widehat A(\kappa)\,$ is referred to as  symbol of $A\,,$
   given as
 $ \widehat A (\kappa ) = \l( A\, \e^{\, -i\kappa x}\r)\,
  \e^{\, +i\kappa x}\,. $}
$\,_xD_\theta^\alpha$ is required to be the logarithm of the
characteristic function of the generic  {\it stable}
(in the L\'evy sense)
probability density, according to the Feller parameterization
\cite{Feller 52,Feller 71}, see also Refs.
\cite{Sato 99,UchaikinZolotarev 99}.
For $\alpha =2$ (hence $\theta=0$) we have
$ \widehat{\,_xD_0^2}(\kappa)  =  -\kappa ^2
= (-i\kappa )^2\,,$ so
we recover the standard second derivative.
More generally for $\theta=0$
we have
$ \widehat{\, _xD_0^\alpha}(\kappa) =
              -|\kappa |^\alpha  = - (\kappa ^2)^{\alpha /2}$ so
$$ \,_xD_0^\alpha  = - \l(-{d^2\over dx^2}\r) ^{\alpha/2}\,.
     \eqno(5) $$
In this case we call  the LHS of Eq.(5)  simply  the
{\it Riesz fractional derivative} operator of order $\alpha \,.$
Assuming $\alpha \ne 1,2$ and taking $\theta$ in its range,
one can show that the explicit expression of the
{\it Riesz-Feller fractional derivative} obtained from Eq. (4)
is
  $$ \,_x D_\theta ^\alpha \, f (x) :=
  -\, \l[
  c_+(\alpha,\theta)\,\,_x  D^{\alpha}_+
+ c_-(\alpha,\theta)\,\,_x  D^{\alpha}_-
     \r]\,  f(x)  \,,  \eqno(6)$$
where
$$
c_+(\alpha,\theta) =
  {\sin \,\l[(\alpha-\theta)\,\pi/2\r] \over\sin\,(\alpha\pi)}\,,
 \qq
   c_-(\alpha,\theta) =
{\sin\,\l[(\alpha+\theta)\,\pi/2\r] \over \sin(\alpha\pi)} \,,
\eqno(7)
$$
and  the
$_xD_\pm^{\alpha} $ are Weyl fractional derivatives
defined as
 $$
 _xD_\pm^{\alpha} \,f(x) = \cases{
     {\ds \pm {d \over dx}} \,\l[\,_x I_\pm^{1-\alpha}\,f(x)\r] \,,
   & if $\q 0<\alpha < 1 \,,$\cr\cr
     {\ds{d^2 \over dx^2}} \,\l[\, _x I_\pm^{2-\alpha}\,f(x)\r]
 \,,
   & if $\q 1<\alpha < 2 \,.$\cr}
\eqno(8)
$$
In Eq. (8) the $\,_x I_\pm^\mu $ ($\mu >0$) denote the Weyl fractional
integrals defined as
 $$ \cases{{\ds _x I_+^\mu  \, f(x)}=
 {\ds \rec{\Gamma(\mu )}}\,
{\ds   \int_{-\infty}^x \!\!  (x-\xi)^{\mu -1}\, f(\xi)\,d\xi} \,,
 \cr\cr
 {\ds _x I_-^\mu  \, f(x)}=
 {\ds \rec{\Gamma(\mu )}}\,
 {\ds  \int_x^{+\infty} \!\! (\xi-x)^{\mu -1}\,f(\xi)\,d\xi} \,.
  \cr} \; (\mu >0) \eqno(9)$$
In the particular case $\theta =0$ we get
$c_+(\alpha ,0)  = c_-(\alpha ,0) =
  {1 /[2 \cos\,(\alpha\pi /2)}] \,,$
and, by passing to the limit for $\alpha  \to 2^-\,,$
we get $\, c_+(2,0) = c_-(2,0) = - 1/2\,. $

For $\alpha =1$ we  have
$$ _xD^1_\theta \, f(x) = \l[ \cos (\theta \pi/2)\,  _xD_0^1
 + \sin (\theta \pi/2)\,  _xD \r] \, f(x) ,\eqno(10)$$
with
$_xD \, f(x) = {\ds {d\over dx}}\, f(x)\,,$ and
$$
  _xD_0^1 \,f(x) = -{d\over dx} \,\l[ _xH\, f(x)\r] \,,
  \q _xH \, f(x) =  {1\over \pi}
\, \l( \int_{-\infty}^{+\infty} {f(\xi )\over x-\xi}\, d\xi\r)
 \,.
\eqno(11) $$
In (11) the operator $\,_xH$ denotes the Hilbert transform
and its singular integral is understood in the Cauchy principal value
sense.

The operator $_xD^\alpha_\theta$
has been referred to as
the {\it Riesz-Feller} fractional derivative since
both Marcel Riesz and William Feller contributed to its definition
\footnote{Originally,   in the late  1940's, Riesz \cite{Riesz 49}
introduced  the pseudo-differential operator $_x I_0^\alpha$ whose
 symbol is $|\kappa|^{-\alpha} \,,$  well defined for any
 positive $\alpha$ with the exclusion of odd integer numbers,
 afterwards named the {\it Riesz potential}.
  The Riesz fractional derivative $_x D_0^\alpha := - \,_x I_0^{-\alpha}$
 defined by analytical continuation was generalized by Feller
 in his 1952 genial paper \cite{Feller 52} to include the skewness
 parameter      of the strictly stable densities.}.

Let us now consider the time-fractional {\it Caputo derivative}.
Following the original idea by Caputo
 \cite{Caputo 67}, see also
\cite{CaputoMaina 71,GorMai CISM97,Podlubny 99},
a proper time fractional derivative of order $\beta \in (m-1,m]$
with $m \in \NN\,, $
useful for physical applications, may be defined in terms
of the following  rule for the Laplace transform
$$ {\cal L} \l\{ _tD_*^\beta \,f(t) ;s\r\} =
      s^\beta \,  \widetilde f(s)
   -\sum_{k=0}^{m-1}    s^{\beta  -1-k}\, f^{(k)}(0^+) \,,
  \q m-1<\beta  \le m \,. \eqno(12)$$
where
$ \widetilde f(s) =
{\cal L} \l\{ f(t);s\r\}
 = \int_0^{\infty} \e^{\ds \, -st}\, f(t)\, dt\,.$
 Then
the {\it Caputo fractional derivative} of
$f(t)$  turns out to be defined as
$$
    _tD_*^\beta \,f(t) :=
\cases{
    {\ds \rec{\Gamma(m-\beta )}}\,{\ds\int_0^t
 {\ds {f^{(m)}(\tau)\, d\tau \over (t-\tau )^{\beta  +1-m}}}} \,,
  & $\; m-1<\beta  <m\,, $\cr\cr
     {\ds {d^m\over dt^m}} f(t)\,,
    & $\; \beta  =m\,. $\cr\cr }
   \eqno(13) $$
In other words the  operator
$\,_tD_*^\beta $ is required to generalize the well-known
rule for the Laplace transform of
the $n$-th order derivative of a given (causal) function keeping the
standard initial value of the function itself
and of its derivatives up to order  $n-1$
\footnote{
The reader should observe that the {\it Caputo} fractional derivative
differs from the usual {\it Riemann-Liouville} fractional derivative
which, defined as the left inverse of the Riemann-Liouville fractional
integral, is  here denoted as
$\, _tD^\beta \,f(t)\,. $
We have,  see \eg \cite{SKM 93},
$$
 _tD^\beta  \,f(t) :=
\cases{
  {\ds {d^m\over dt^m}}\,\l[
  {\ds \rec{\Gamma(m-\beta )}\,\int_0^t
    {f(\tau)\,d\tau  \over (t-\tau )^{\beta  +1-m}} }\r] \,,
 & $\; m-1 \le \beta  < m\,,$ \cr\cr
     {\ds {d^m\over dt^m}} f(t)\,,
    & $\; \beta  = m\,. $\cr\cr }
 $$
When $\beta $ is not integer and both fractional derivatives exist
we have between them the following relation,
 see \eg \cite{GorMai CISM97},
$$ _tD_*^\beta  \,f(t)  \, = \, _tD^\beta  \,\l[ f(t) -
  \sum_{k=0}^{m-1} f^{(k)}(0^+)\,{t^k\over k!} \r]\,,\q
    m-1 <\beta <m\,,$$
or
$$ _tD_*^\beta  \,f(t)  \, = \, _tD^\beta  \, f(t) -
    \sum_{k=0}^{m-1}    f^{(k)}(0^+) \,
{t^{k-\beta }\over \Gamma(k-\beta+1)}\,,
\q       m-1 <\beta <m\,. $$
The {\it Caputo} fractional derivative,
 originally ignored in the  mathematical treatises,
represents a sort of regularization in the time origin for the
{\it Riemann-Liouville} fractional derivative
and  satisfies the  relevant property
of being zero when applied to a constant.
For more details on this fractional derivative
we refer the interested reader to
Gorenflo and Mainardi \cite{GorMai CISM97} and Podlubny \cite{Podlubny 99}.
\\ {\it Added Note}: Nowadays the reader can find an exhaustive treatment of fractional derivatives
in  the treatise by A.A. Kilbas, H.M. Srivastava  and J.J. Trujillo,
{\it Theory and Applications of Fractional Differential Equations},
Elsevier, Amsterdam (2006).}.

The fundamental solution (or the {\it Green function})
$\G(x,t)\,$
of Eq. (2)
is the solution corresponding to the initial condition
$u(x,0^+) = \delta (x)\,.$
We note that if $1 < \beta \le 2$ we  add   the condition
$u_t(x,0^+)=  0\,. $

In the particular case of the standard diffusion equation (1)
the Green function is nothing but
the Gaussian probability density function with variance
$\sigma^2  =2t\,,$ namely
$$G_{2,1}^0 (x,t)
 = {1\over 2\sqrt{\pi }}\,t^{-1/2}\, \e^{-\ds x^2/(4t)}\,.
\eqno(14)$$
We note that in the limiting case  $\{\alpha =\beta =2\}$
we recover the D'Alembert {\it wave equation} with
Green function   $G_{2,2}^0 (x,t) = [\delta(x+t) + \delta (x-t)]/2\,. $

In the general case the application of the transforms of Fourier
and Laplace  in succession
to Eq (2) yields
 $$- \psi_\alpha^\theta(\kappa ) \,\Gks  \,= \,
 s^\beta \,\Gks - s^{\beta -1}    \,,
\eqno(15)  $$
namely
$$
  \Gks = {s^{\beta -1} \over
 s^\beta + \psi_\alpha^\theta (\kappa)} \,. \eqno(16)$$
By using the known scaling rules
for the Fourier  and Laplace  transforms,
following the arguments by Mainardi {\it et al.} 
\cite{Mainardi LUMAPA01}
one can prove without inverting Eq. (16) that $\G(x,t)$
has the  scaling property
$$ \G(x,t)  =
    t^{-\beta /\alpha}\,\K \l(x/t^{\beta/\alpha}\r)\,.
   \eqno(17)     $$
Here $ x/t^{\beta/\alpha}\,$ acts as the similarity variable  and
$\K (\cdot)$ as  the {\it reduced Green function}.
For the analytical and computational determination of
the  reduced Green function
we can restrict the attention to $x>0\,$
because of the
{\it symmetry relation} 
$ K_{\alpha ,\beta}^\theta(-x) = K_{\alpha ,\beta}^{-\theta}(x)
\,.$
Extending the paper by Gorenflo {\it et al.} 
\cite{GoIsLu 00} where the space-time fractional diffusion equation
has been considered with the restriction $\{1<\alpha \le 2,\,
\theta =0,\, 0<\beta \le 2\}$,
Mainardi {\it et al.} 
\cite{Mainardi LUMAPA01}
have first inverted the Laplace transform and then
have obtained the following
 Fourier integral representation of the
reduced Green function
$$ \K (x) =      \rec{2\pi}\,
\int_{-\infty}^{+\infty} \e^{\,\ds -i \kappa x} \,
  E_\beta \l[- \psi_\alpha ^\theta(\kappa ) \r] \, d\kappa
 \,.
\eqno(18) $$
Here
$E_{\beta}$ denotes the entire transcendental function, known as
the Mittag-Leffler function of order $\beta\,,$
defined in the complex plane by the power series
$$ E_\beta (z) :=
    \sum_{n=0}^{\infty}\,
   {z^{n}\over\Gamma(\beta\,n+1)}\,, \q \beta >0\,, \q z \in \CC\,.
 \eqno  (19)$$
For detailed information on the Mittag-Leffler function
the interested reader may  consult \eg
\cite{Erdelyi HTF} (Vol. 3, Ch. 18, pp. 206-227)
and \cite{GorMai CISM97,Podlubny 99,SKM 93,MaiGor JCAM00}.
By using the convolution theorem for the Mellin transform
\footnote{If
$$
   {\cal M} \, \{ f(r ); s\} = f^*(s)=
   \int_0^{+\infty} f(r)\,
 r^{s-1}\,  dr,  \q  \gamma_1< \Re\, (s) <\gamma_2
$$
denotes the Mellin transform of 
$f(r)\,,$ the inversion is provided by
$$
 {\cal M}^{-1}\, \{  f^*(s ); r \} =f(r)=
{1\over 2\pi i}\int_{\gamma -i
\infty}^{\gamma +i\infty} f^*(s)\, r^{-s} \,ds
$$
where $\ r>0\,,$ $\, \gamma = \Re\,(s) \,,$
$\, \gamma_1< \gamma <\gamma_2\,.$
The Mellin convolution formula reads
$$ h(r ) = \int\limits_0^\infty
 \rec{\rho }\, f(\rho)\,g(r/\rho  )\,{d\rho  }
\,\stackrel{{\cal M}}{\leftrightarrow}\,
 h^*(s) = f^*(s)\,g^*(s) \,.       $$},
the authors in Ref. \cite{Mainardi LUMAPA01}  have provided the Mellin-Barnes
integral representation
\footnote{
The names refer to the two authors,
who in the 
beginning of the past century
developed the theory of these integrals  using them
for a complete integration of the hypergeometric differential equation.
However, as pointed out  in \cite{Erdelyi HTF}
(Vol. 1, Ch. 1, \S 1.19, p. 49), these integrals were first used
by S. Pincherle in 1888. For a revisited analysis of the pioneering work
of Pincherle (Professor of Mathematics at the
University of Bologna from 1880 to 1928) we refer
to the recent paper by
Mainardi and Pagnini \cite{MainardiPagnini OPSFA01}.
As a matter of fact this type of integrals turns out to be useful
in inverting the Mellin transform.}
$\;$ for the general case $0<\alpha \le 2$, $0<\beta\le 2$
as in Eqs. (2) and (3):
$$  \K(x) =
{1\over  \alpha x}
{1\over 2\pi i} \int_{\gamma-i\infty}^{\gamma+i\infty}
{\Gamma({s\over \alpha}) \, \Gamma(1-{s\over \alpha}) \,\Gamma(1-s)
 \over \Gamma(1-{\beta\over \alpha}s) \,
 \Gamma ( \rho \,s)\,
 \Gamma (1-\rho \,s)}
 \, x^{\,\ds s}\,  ds,
\; \rho =   { \alpha -\theta \over 2\,\alpha },
\eqno(20)  $$
where $0< \gamma < \hbox{min} \{\alpha ,1\}.$
Then, 
 after changing $s$
into $-s$ and using the known property $\Gamma(1+z) = z\,\Gamma(z)\,,$
we can  write the Mellin transform of
$x \, K_{\alpha ,\beta}^\theta(x)$ as
$$  \int_0^{+\infty}  \!\!
 K_{\alpha ,\beta}^\theta(x) \,x^{\, \ds s} \,dx =
 \rho \,
 { \Gamma(1-{s/\alpha})\,\Gamma(1+{s/\alpha}) \,\Gamma(1+s)
 \over
 \Gamma (1-\rho \,s)\, \Gamma (1+\rho \,s)\,\Gamma(1+{\beta\,s/ \alpha})}
\,,
 \eqno(21) $$
where
$-\hbox{min}\{\alpha, 1\}< \Re (s) <\alpha \,. $
In particular 
we find
$  \int_0^{+\infty}   K_{\alpha ,\beta}^\theta(x)  \,dx = \rho \,$
(with $ \rho  = 1/2$ if $\theta =0$).
We note that Eq. (21) is strictly valid as soon as cancellations in the
"gamma fraction" at the RHS are not possible.
Then this equation  allows us to evaluate (in $\RR_0^+$)
the (absolute) moments of order $\delta  $ for the Green function
if   $ -\hbox{min}\{\alpha, 1\} <\delta  <\alpha\,.  $
In other words, it states that
$\,K_{\alpha ,\beta}^\theta(x)= {\cal O} \l(x^{-(\alpha +1)}\r)$
 as $x \to +\infty\,.$
When cancellations occur  in the "gamma fraction"  the range
of $\delta$ may change. 
An interesting case is 
$\,\{ \alpha =2, \, \theta =0, \, 0<\beta < 2\}$
 ({\it time-fractional diffusion} including {\it standard diffusion}),
where Eq. (21)
reduces to $$ \int_0^{+\infty}  \!\!
 K_{2,\beta}^0(x) \,x^{\, \ds s} \,dx =
 \rec{2} \,
 {\Gamma(1+s)
 \over
\Gamma(1+{\beta\,s/2})}
\,,\q \Re (s) >-1\,.
 \eqno(22) $$
This result  proves
the existence of all moments of order $\delta  >-1$
for the corresponding Green function which reads
$$K_{2,\beta}^0(x) =
 {1\over 2} M_{\beta/2} (x) =
{1\over  2 x}
{1\over 2\pi i} \int_{\gamma-i\infty}^{\gamma+i\infty}
{\Gamma(1-s)
 \over   \Gamma (1- \beta s/2)}
 \, x^{\,\ds s}\,  ds\,,\q x>0\,,\eqno(23)  $$
where $M_{\beta/2}$ denotes  a function
of the Wright type  introduced by Mainardi
\cite{Mainardi CHAOS96,Mainardi CISM97},
see also \cite {Podlubny 99,GoLuMa 99,GoLuMa 00}.
For the case $\alpha =\beta $ (that we call
{\it neutral-fractional diffusion})
we obtain from (20) an elementary  expression
$$ \qq \qq K_{\alpha,\alpha}^\theta (x) = N_\alpha ^\theta(x) =
{1\over  \alpha x}
{1\over 2\pi i} \int_{\gamma-i\infty}^{\gamma+i\infty}
{\Gamma({s\over \alpha}) \, \Gamma(1-{s\over \alpha})
 \over
 \Gamma ( \rho \,s)\,
 \Gamma (1-\rho \,s)}
 \, x^{\,\ds s}\,  ds  \qq \qq
\eqno(24)$$
$$  =
{1\over  \alpha x}
{1\over 2\pi i} \int_{\gamma-i\infty}^{\gamma+i\infty}
 {\sin (\pi\, \rho \,s)
 \over
 \sin (\pi \,s/\alpha )}
 \, x^{\,\ds s}\,  ds
=     {{1\over\pi}}\,{x^{\alpha-1} \sin[{\pi\over 2}(\alpha -\theta )] \over
1 + 2x^\alpha \cos[{\pi\over 2}(\alpha -\theta)] + x^{2\alpha}}
\,.
     $$
where $0<\gamma <\alpha \,. $

 We note that the
 Mellin-Barnes integral representation allows us to
 construct computationally   the fundamentals solutions
of Eq. (2) for any triplet $\{\alpha ,\beta, \theta\}$
($\alpha \ne \beta $)
by  matching their
convergent and asymptotic expansions, as shown in
\cite{Mainardi LUMAPA01}. 
Readers acquainted with  Fox $H$ functions can recognize
in Eq. (20) the representation of a certain function of this
class, see \eg
 \cite{UchaikinZolotarev 99, SKM 93,MathaiSaxena H,Srivastava H}.
Unfortunately, as far as we know, computing routines for this general class
of special functions are not yet available\footnote{
{\it Added Note} Nowadays the reader can find the representation of the fundamental solutions
in terms of Fox $H$ functions in the paper by F. Mainardi, G. Pagnini and R.K. Saxena,
 Fox $H$ functions in fractional diffusion,
 {\it   J. Computational and Applied  Mathematics} {\bf 178},  321-331 (2005).}.

Let us now point out that, in the peculiar cases
of {\it space-fractional-diffusion},
{\it time-fractional-diffusion},
and {\it neutral-fractional-diffusion},
non-negativity of the corresponding Green functions
can be proven; being
normalized in $\RR$  these functions can indeed be interpreted as
probability densities.
In particular,
for  $\beta =1$ and $0<\alpha <2$
({\it strictly space-fractional diffusion})
we recover
the class  of the strictly stable (non-Gaussian)
densities $\,L_\alpha ^\theta(x)\,$
exhibiting fat tails (with the algebraic decay
$\propto |x|^{-(\alpha +1)}$)
  and infinite variance,
whereas for $\alpha =2$ and $0<\beta <1$
({\it strictly time-fractional diffusion})
the class
of the Wright-type densities $\, M_{\beta/2}(x)/2\,$
exhibiting stretched exponential
tails   and finite variance proportional to $t^\beta \,.$


The meaning of probability density
can  be extended under proper conditions  to the reduced Green function
of the  {\it space-time-fractional} partial differential equation
(2) in virtue of the    following identity, proven
in Ref. \cite{Mainardi LUMAPA01},
$$ \K(x) = \cases{
   \alpha   \, {\ds \int_0^\infty} \!\!
  \l[\xi ^{\alpha -1}\,  {M}_{\beta}\l(\xi ^{\alpha}\r)\r]\,
  {L}_{\alpha}^\theta\l({x/\xi }\r) \,
      {\ds{d\xi  \over \xi }} \,,&$\; 0<\beta <1\,,$ \cr\cr
  {\ds \int_0^{\infty}} \!\!
M_{\beta / \alpha}(\xi ) \,N_\alpha^\theta(x/\xi)\,
    {\ds{d\xi  \over \xi} }\,, & $\; 0<\beta/\alpha<1 \,.$ \cr}
\eqno(25)$$
Then, due to the previous discussion,
in the cases
$\{0<\alpha <2\,,\, 0<\beta < 1 \}\,$
{and} $\{1<\beta \le \alpha < 2\} $
we obtain  a class of  probability densities
(symmetric or non-symmetric according to $\theta =0$ or $\theta \ne 0$)
which exhibit  fat tails
with an algebraic decay $\propto |x|^{-(\alpha +1)}\,.$
Thus, they belong to the domain of attraction of the L\'evy stable densities
of index $\alpha $ and  can be referred to as
{\it fractional stable densities}.

\section{CONCLUSIONS}

The above analysis of the Cauchy problem for the
{\it space-time fractional} diffusion equation  shows the fundamental
and combined roles of the integral transforms of Fourier, Laplace and
Mellin type. By the aid of the transforms of Fourier and Laplace
 the
scaling and similarity properties of the Green function can be
easily derived.
The Mellin transform
allows us  to obtain for the Green function
a general  representation  in terms
of Mellin-Barnes integrals (hence a computational form
in terms of convergent and asymptotic series),
 and the extension of its
probability interpretation\footnote{
{\it Added Note} It is worthwhile to attract  the reader's attention
to the role of the Mellin transform to derive 
subordination laws in fractional diffusion  processes, see
F. Mainardi, G. Pagnini and R. Gorenflo,
Mellin transform and subordination laws in fractional diffusion  processes,
   {\it Fractional Calculus and Applied Analysis} {\bf 6} No. 4,  441-459 (2003).
 [E-print {\tt http://arxiv.org/abs/math/0702133}]}


\vspace*{-2pt}

\section*{Acknowledgements}

The author is very grateful to Rudolf  Gorenflo
for  inspiring discussions and helpful comments.
He is also grateful to the National Group of Mathematical Physics
(G.N.F.M. - I.N.D.A.M.) and the National Institute of Nuclear Physics
(I.N.F.N. - Sezione  di Bologna) 
for partial support.
This article is based on an invited lecture at the 3rd International
ISAAC Congress, Freie Universit\"at Berlin, 20-25 August 2001
(Sub-session 1.3: Integral Transforms and Applications).


\end{document}